\documentclass[12pt,a4paper]{amsart}
\usepackage{graphicx,amssymb}
\input xy
\xyoption{all}
\usepackage[all]{xy}
\usepackage{hyperref}

\newtheorem{thm}{Theorem}[section]
\newtheorem{cor}[thm]{Corollary}
\newtheorem{lem}[thm]{Lemma}
\newtheorem{prop}[thm]{Proposition}
\newtheorem{hyp}[thm]{Hypothesis}
\theoremstyle{definition}

\newtheorem{rem}[thm]{Remark}

\numberwithin{equation}{section}
\newcommand{\QQ}{\mathbb Q}

\newcommand{\ZZ}{\mathbb Z}
\newcommand{\CC}{\mathbb C}
\newcommand{\PP}{\mathbb P}

\newcommand{\lra}{\longrightarrow}

\newcommand{\ra}{\rightarrow}

\newcommand{\cW}{\mathcal{W}}
\newcommand{\cV}{\mathcal{V}}
\newcommand{\cD}{\mathcal{D}}
\newcommand{\cT}{\mathcal{T}}

\newcommand{\cK}{\mathcal{K}}

 \DeclareMathOperator{\Gal}{Gal}

 \DeclareMathOperator{\Nm}{{Nm}}
  \DeclareMathOperator{\Prym}{{Prym}}

 \DeclareMathOperator{\Stab}{Stab}
 \DeclareMathOperator{\tr}{tr}

 \DeclareMathOperator{\im}{Im}
 \DeclareMathOperator{\Div}{{Div}}

\DeclareMathOperator{\End}{{End}}

\DeclareMathOperator{\Alt}{{Alt}}

\begin{document}

\title[ ]{Prym-Tyurin varieties via Hecke algebras}
\author{Angel Carocca, Herbert Lange, Rub\'i E. Rodr\'iguez and Anita M. Rojas}

\address{A. Carocca\\Facultad de Matem\'aticas, Univ. Catolica, Casilla 306-22, Santiago, Chile}

\email{acarocca@mat.puc.cl}

\address{H. Lange\\Mathematisches Institut,
              Universit\"at Erlangen-N\"urnberg\\Germany}

\email{lange@mi.uni-erlangen.de}

\address{R. E. Rodr\'iguez\\Facultad de Matem\'aticas, Univ. Catolica, Casilla 306-22, Santiago, Chile}

\email{rubi@mat.puc.cl}

\address{A. M. Rojas \\Departamento de Matematicas, Facultad de Ciencias,
Universidad de Chile, Santiago\\Chile}
\email{anirojas@uchile.cl }

\thanks{The first, third and fourth author were supported by Fondecyt grants 1060743, 1060742 and 11060468 respectively}%

%\date{ }%
%\dedicatory{ }%
%\commby{ }%
% ----------------------------------------------------------------
\begin{abstract}
Let $G$ denote a finite group and $\pi: Z \ra Y$ a Galois covering
of smooth projective curves with Galois group $G$. For every
subgroup $H$ of $G$ there is a canonical action of the
corresponding Hecke algebra $\QQ[H \backslash G/H]$ on the
Jacobian of the curve $X = Z/H$. To each rational irreducible
representation $\mathcal{W}$ of $G$  we associate an idempotent in
the Hecke algebra, which induces a correspondence of the curve $X$
and thus an abelian subvariety $P$ of the Jacobian $JX$. We give
sufficient conditions on $\mathcal{W}$, $H$, and the action of $G$
on $Z$ for $P$ to be a Prym-Tyurin variety. We obtain
many new families of Prym-Tyurin varieties of arbitrary exponent
in this way.
\end{abstract}

\maketitle

%\tableofcontents

\section{Introduction}

A {\it Prym-Tyurin variety of exponent $q$} is by definition a principally polarized abelian variety $(P,\Xi)$, for which
there exists a smooth projective curve $C$ and an embedding $P \hookrightarrow JC$ into the Jacobian of $C$ such that the
restriction of the canonical polarization $\Theta$ of $JC$ is the $q$-fold of $\Xi$:
$$
\Theta|P = q\Xi.
$$
One point of interest in these varieties is that the structure of a Prym-Tyurin variety allows to study the geometric properties of the underlying abelian varieties via curve theory.

According to the Theorem of Matsusaka-Ran \cite[11.8.1]{bl} Prym-Tyurin varieties of exponent 1 are exactly the
canonically polarized products of Jacobians. Welters showed in \cite{w} that Prym-Tyurin varieties
of exponent 2 are exactly the classical
Prym varieties or some of their specializations. The Abel-Prym-Tyurin map of $(P,\Xi)$ is the following composition of maps
$$
\alpha_P: C \stackrel{\alpha_C}{\lra} JC = {\widehat {JC}} \ra \widehat P = P,
$$
where $\alpha_C$ denotes the usual Abel map and $\widehat A$ the dual of an abelian variety $A$. The property of being a
Prym-Tyurin variety can then be expressed by certain properties of the curve $\alpha_P(C)$ in $P$. One can deduce
a criterion (see \cite[Welters' criterion 12.2.2]{bl}) for $(P,\Xi)$ to be a Prym-Tyurin variety in term of a curve in $P$. Using
this it is easy to see, using complete intersection curves on the Kummer variety of $P$, that every principally polarized
abelian variety of dimension $g$ is a Prym-Tyurin variety of exponent $2^{g-1}(g-1)!$.

The problem is to construct Prym-Tyurin varieties of small exponent $\geq 3$. There are in fact not many examples. A series
of examples was given by Kanev in \cite{k} using Weyl groups of type $A_n, D_n, E_6$ and $E_7$. Generalizing the
construction of \cite{lrr}, Salomon gave in \cite{so} examples of Prym-Tyurin varieties of arbitrary exponent using
some special graph constructions. There are
a few other examples which shall not be mentioned here.

The construction of these examples are based on the following two facts: (i) The use of correspondences: The set of
abelian subvarieties of a polarized abelian variety is in canonical bijection to the set of symmetric (with respect
to the Rosati involution) idempotents of its endomorphism algebra. On the other hand, it follows almost from the definitions
that the obvious map
$$
\Div_{\QQ}(C \times C) \lra \End_{\QQ}(JC)
$$
from the algebra of rational correspondences of $C$ to the endomorphism algebra of its Jacobian is surjective. Hence one can
describe an abelian subvariety of $JC$ by correspondences. (ii) Kanev's Criterion: In \cite{k1} Kanev showed that an abelian
subvariety of $JC$ is a Prym-Tyurin variety of exponent $q$ if it is given by an integral effective fixed-point free correspondence
on $C$ whose associated endomorphism of $JC$ satisfies the equation
\begin{equation} \label{eq1.1}
x^2 + (q-2)x - (q-1) = 0.
\end{equation}
This is a direct generalization of Wirtinger-Mumford's construction of classical Prym varieties and, in fact,
most Prym-Tyurin varieties are given by constructing integral symmetric fixed point free correspondences satisfying \eqref{eq1.1}.
There is an analogous result for correspondences with some fixed points, due to Ortega \cite{o}, which also can be applied.\\

Our construction of such correspondences is, roughly speaking, as follows: Let $G$ denote a finite group and $\pi: Z \ra Y$ be
a Galois covering of smooth projective curves with Galois group $G$. The group action induces a homomorphism of the group algebra $\QQ[G]$
into the endomorphism algebra $\End_{\QQ}(JZ)$. We identify the elements of $\QQ[G]$ with their images in $\End_{\QQ}(JZ)$.
Hence any $\alpha \in \QQ[G]$ defines an abelian subvariety $\im(\alpha)$ of $JZ$. For the details we refer to Section 3.1.
On the other hand, every element $\alpha = \sum_{g \in G} \alpha_g g \in \QQ[G]$ defines a rational correspondence on $Z$, namely
$$
\cD_{\alpha} = \sum_{g \in G} \alpha_g \Gamma_g \in \Div_{\QQ}(Z \times Z).
$$
where $\Gamma_g$ denotes the graph of the automorphism $g$ of $Z$. The problem is to find integral symmetric fixed-point
free correspondences satisfying \eqref{eq1.1} in this way. For this we proceed as follows:

For simplicity we first describe
a special case. Let $V$ be an irreducible complex representation of $G$ and $\cW$ the associated rational irreducible representation.
For any subgroup $H$ of $G$, the covering $\pi$ factorizes via a covering $\varphi_H: X \lra Y$ where $X$ denotes the quotient $Z/H$.
The Hecke algebra $\QQ[H \backslash G/H]$ is a subalgebra of $\QQ[G]$ and the above homomorphism induces a homomorphism
$$
\QQ[H \backslash G/H] \lra \End_{\QQ}(JX).
$$
which factorizes via $\Div_{\QQ}(X \times X)$. Then we associate to the pair $(H,\cW)$ an integral symmetric correspondence
${\overline \cD}_{H,\cW}$ on $X$ in a canonical way. It is defined via the projectors of $\QQ[G]$ given in
\cite{caro}. Our essential assumption is
\begin{equation} \label{eq1.2}
 H \; \mbox{is a subgroup satisfying} \; \dim V^H = 1 \; \mbox{and maximal with this property},
\end{equation}
where $V^H$ denotes the fixed subspace of $V$ under the action of $H$.

In fact, under this assumption we can show (Proposition \ref{prop3.9}) that a modified version $\cK_X$ of
${\overline \cD}_{H,\cW}$ is an effective integral
symmetric correspondence on $X$, which in the special case $Y = \PP^1$ satisfies equation \eqref{eq1.1}. Moreover we can
compute its fixed points in terms of the ramification of $\pi: Z \ra \PP^1$.

In order to formulate the result, we use
the following notation: We denote by $K_V$ the field generated by the values of the character of $V$. Moreover,
if $C$ denotes a conjugacy class of cyclic subgroups of the group $G$, a branch point $y \in Y$ of the covering $\pi$ is called
{\it of type} $C$,
if the stabilizer of any point $z$ in the fibre $\pi^{-1}(y)$ is a subgroup in the class $C$. We define {\it the geometric signature} of the covering
$\pi: Z \ra Y$ to be the tuple $[g;(C_1,m_1), \ldots, (C_t,m_t)]$, where $g$ is the genus of the quotient curve $Y$,
the covering $\pi$ has a total of $\sum_{j=1}^t m_j$ branch points, and exactly $m_j$ of them are of type $C_j$
for $j = 1, \ldots, t$ (see \cite{R}). Using this we prove the following\\

\noindent
{\bf Theorem.} Let {\it $\cW$
denote a nontrivial rational irreducible representation of $G$,
with associated complex irreducible representation $V$, and $H$ a subgroup of $G$ satisfying \eqref{eq1.2}.
Suppose that the action of $G$ has geometric signature $[0;(C_1,m_1), \ldots , (C_t,m_t)]$
satisfying
\begin{equation}  \label{eqn1.3}
\sum_{j=1}^t m_j \left( q[K_V:\QQ](\dim V - \dim V^{G_j})  -
([G:H] -|H \backslash G/G_j|) \right) = 0,
\end{equation}
for an integer $q$ given in terms of the algebraic data, where $G_j$ is of class $C_j$,
Then the correspondence $\cK_X$ defines a Prym-Tyurin variety of exponent $q$. }\\

Using the theorem one obtains Kanev's examples and several more,
some of which are included in Sections 5.1 and 5.2. However we
could not find examples that exhibited the full force of the
theorem; we were interested in finding Prym-Tyurin varieties of
small exponent bigger than $2$, using complex representations
whose field of definition properly contains $\QQ$, as these would
provide entirely new examples, different from the existing ones.
In fact, we wrote a computer program which for many groups $G$ and
all of their subgroups and irreducible complex representations
such that $\QQ\subsetneq K_V$ verifies the hypothesis of the
theorem. In this way we found many examples of exponents $1$ and
$2$, which however are not interesting from our point of view.

The new idea was to start, instead of one irreducible representation, with several representations satisfying some additional
conditions (see Hypothesis \ref{hyp}) and generalize the above result to Theorem \ref{thm4.8}. This gave in fact many examples, most of them new,
and which also include the above mentioned examples of Salomon. For details see Section 5.3 and the forthcoming paper \cite{clrr},
in which more examples are given and these Prym-Tyurin varieties will be investigated in detail.\\

The contents of the paper are as follows. Section \ref{S:2} contains some algebraic preliminaries. The essential result is
Proposition \ref{prop2.4}, which describes the coefficients of a certain idempotent of the group algebra in term of the character of $V$.
In Section \ref{S:3} we define the correspondence $\cK_X$ and derive its properties. Section \ref{S:4} contains the proof of our main result Theorem
\ref{thm4.8}. We also compare our construction with the construction of \cite{k}. Finally in Section \ref{S:5} we give some examples.

We suppose throughout that the curves are defined over any algebraically closed field of characteristic 0.
Moreover the curves always will be smooth and projective.

\section{Algebraic preliminaries}\label{S:2}

\subsection{The group algebra}
Let $G$ be a finite group. In order to fix the notation, we start by recalling some basic properties of representations
of $G$ (see \cite{cr}). For any field $K$ of characteristic 0 we denote by $K[G]$ the group algebra of $G$ over $K$.
It is a semisimple algebra, whose simple components correspond
one-to-one to the irreducible $K$-representations of $G$. We may identify the elements of $K[G]$ with the $K$-valued
functions on $G$. From this point of view, the multiplication in $K[G]$ is convolution
$$
(f_1f_2)(g) = \sum_{g_1g_2 = g} f_1(g_1) f_2(g_2)
$$
for $f_1, f_2 \in K[G]$ and $g \in G$.
In this paper the field $K$ will be either $\CC$ or $\QQ$.
% A scalar product on $\CC[G]$ is defined by
%$$
%\langle f_1,f_2 \rangle = \frac{1}{|G|} \sum_{g \in G} f_1(g) \overline{f_2(g)},
%$$

For any complex irreducible representation $V$ of $G$, we denote by $\chi_V$ its
character, by $L_V$ its field of definition and by $K_V$ the subfield $K_V = \QQ(\chi_V(g)\;|\;g \in G)$. $L_V$ and $K_V$ are
finite abelian extensions of $\QQ$. We denote by
$$
m_V = [L_V:K_V]
$$
the Schur index of $V$. For any automorphism $\varphi$ of $L_V/\QQ$ we denote by $V^{\varphi}$ the
representation conjugate to $V$ by $\varphi$.

If $\cW$ is a rational irreducible representation of $G$,
then there exists a complex irreducible representation $V$ of $G$, uniquely determined up to conjugacy in $\Gal(L_V/\QQ)$,
such that
$$
\cW \otimes_{\QQ} \CC  \simeq \bigoplus_{\varphi
\in \Gal(L_V/\QQ)} V^{\varphi} \simeq m_V \bigoplus_{\tau \in
\Gal(K_V/\mathbb{Q})} V^{\tau}
$$
We call $V$ the complex irreducible representation {\it associated to}
$\mathcal{W}$.

The central idempotent $e_V$ of $\CC[G]$ that generates the simple subalgebra of $\CC[G]$ corresponding to $V$ and
the central idempotent $e_{\cW}$ of $\QQ[G]$ that generates the simple subalgebra of $\QQ[G]$ corresponding to $\cW$ are given by
\begin{equation} \label{eq2.1}
e_V = \frac{\dim V}{|G|} \sum_{g \in G} \chi_V(g^{-1})g \quad \mbox{and} \quad e_{\cW} =
\frac{\dim V}{|G|} \sum_{g \in G} \tr_{K_V/\QQ} (\chi_V(g^{-1}))g.
\end{equation}

\begin{lem} \label{lem2.1}
If $e$ is an idempotent of $K[G]$, the map
$$
\End_{K[G]}(K[G]e) \ra eK[G]e,\quad \varphi \mapsto \varphi(e)
$$
is an anti-isomorphism of $K$-algebras.
\end{lem}

\begin{proof}
Note first that $\varphi(e) = \varphi(e\cdot e) = e \cdot \varphi (e) \in eK[G]e$. For the proof that the map is ``anti'',
suppose that $\varphi_i(e) = f_ie$ for $i = 1, 2$. Then
$$
\varphi_1 \varphi_2(e) = \varphi_1 (\varphi_2(e)) = \varphi_1 (f_2e) =  f_2 \varphi_1 (e) = f_2ef_1e=
\varphi_2(e) \varphi_1(e).
$$
Finally, for any $efe \in eK[G]e$, the map $\varphi \in \End_{K[G]}(K[G]e)$,
defined by $\varphi(f'e) = f'efe$ for all $f' \in K[G]$,
maps to $efe$. Hence the map is bijective.
\end{proof}

\subsection{The Hecke algebra of a subgroup} \label{ss2.2}
Let $H$ be a subgroup of $G$.
The element
\begin{equation} \label{eq2.3}
p_H = \frac{1}{|H|} \sum_{h \in H} h
\end{equation}
is the central idempotent of $K[H]$ corresponding to the trivial representation of $H$. Moreover, the left ideal
$K[G]p_H$ defines the $K$-representation of $G$ induced by the trivial representation of $H$ for any field $K$ as above.
In the sequel we
denote this representation by $\rho_H$.

The $K$-algebra $p_HK[G]p_H$, considered as a subalgebra of $K[G]$, consists of the $K$-valued functions on $G$ which
are constant on each double coset $HgH$ of $H$ in $G$. It is called the {\it Hecke algebra} over $K$ of $H$ in $G$, and it is
usually denoted by $K[H \backslash G/H]$.

The idempotent $e_V$ of \eqref{eq2.1} is central in $\CC[G]$. This implies that the element
\begin{equation} \label{eq2.4}
f_{H,V} = p_H e_V = e_V p_H
\end{equation}
is an idempotent of $\CC[G]$ (or zero) which satisfies
\begin{itemize}
\item $hf_{H,V} = f_{H,V} = f_{H,V}h$ for all $h \in H$,
\item the left ideal $\CC[G]f_{H,V}$ defines the representation $V$ with multiplicity $\dim V^H$.
\end{itemize}
The last property follows from the fact that $\rho_H \simeq \CC[G]p_H$ and $V$ occurs in $\rho_H$ with multiplicity $\dim V^H$.\\

Similarly, since the idempotent $e_{\cW}$ of \eqref{eq2.1} is central in $\QQ[G]$, the element
\begin{equation} \label{eq2.5}
f_{H,\cW} := p_H e_{\cW} = e_{\cW} p_H
\end{equation}
is an idempotent of $\QQ[G]$ (or zero) which satisfies
\begin{itemize}
\item $hf_{H,\cW} = f_{H,\cW} = f_{H,\cW} h$ for all $h \in H$,
\item the left ideal $\QQ[G]f_{H,\cW}$ defines the representation $\cW$ with multiplicity $\frac{\dim V^H}{m_V}$.
\end{itemize}
The last property and the fact that $\frac{\dim V^H}{m_V}$ is an integer follow from the equation
$\dim V^H = \langle \rho_H,V \rangle$ and the
fact that $\rho_H$ is a rational representation.  For a complete proof see \cite[Theorem 4.4]{caro}.\\

According to \cite[Theorem 5.1.7]{h} there exists a set of representatives
$$
\{ g_{ij} \in G \;| \; i = 1, \ldots, d \; \mbox{and} \; j = 1, \ldots, n_i \}
$$
for both the left cosets and right cosets of $H$ in $G$, such that
$$
G = \sqcup_{i=1}^d Hg_{i1}H \quad \mbox{and} \quad Hg_{i1}H = \sqcup_{j=1}^{n_i} g_{ij}H = \sqcup_{j=1}^{n_i} Hg_{ij}
$$
are the decompositions of G into double cosets, and of the double cosets into right and left cosets of $H$ in $G$. Moreover,
we assume $g_{11} = 1_G$.

Clearly $\chi_V(g^{-1}) = \overline{\chi_V(g)}$ is an algebraic integer for all $g \in G$ (\cite[Corollary 30.11]{cr}). This
implies that
$\tr_{K_V/\QQ}(\chi_V(g)) = \tr_{K_V/\QQ}(\chi_V(g^{-1}))$ is an integer for every $g \in G$. We conclude
\begin{equation}  \label{eq2.2}
\sum_{h \in H} \tr_{K_V/\QQ}(\chi_V(hg_{ij}^{-1})) = \sum_{h \in H} \tr_{K_V/\QQ}(\chi_V(g_{ij}h)),
\end{equation}
and use this to prove the following lemma.

\begin{lem} \label{lem2.2} For any $1 \leq i \leq d$ and $1 \leq j \leq n_i$ we have
\begin{equation} \label{eq2.6}
a_i :=\sum_{h \in H} \tr_{K_V/\QQ}(\chi_V(hg_{ij}^{-1})) = \sum_{h \in H} \tr_{K_V/\QQ}(\chi_V(hg_{i1}^{-1})) \in \ZZ.
\end{equation}
\end{lem}
\begin{proof}
Since $g_{ij}$ and $g_{i1}$ are in the same double coset of $H$ in $G$, there exist $k, k' \in H$ such that $g_{ij} = kg_{i1}k'$.
Hence for all $h \in H$ we have $g_{ij}h = kg_{i1}(k'hk)k^{-1}$, which implies
$$
\chi_V(g_{ij}h) = \chi_V(g_{i1}(k'hk))
$$
for all $h \in H$. Now the map $g_{ij}H \ra g_{i1}H, \; g_{ij}h \mapsto g_{i1}k'hk$ is a bijection, which gives
\begin{equation} \label{eqn2.3}
\sum_{h \in H} \chi_V(g_{ij}h) = \sum_{h \in H} \chi_V(g_{i1}h).
\end{equation}
Now applying $\tr_{K_V/\QQ}$ and \eqref{eq2.2} twice gives the assertion.
\end{proof}

Now consider the following element of $\frac{1}{|H|} \ZZ[G]$
$$
F_{H,\cW} := \sum_{i=1}^d a_i \sum_{j=1}^{n_i} g_{ij} p_H.
$$

\begin{lem}   \label{lem2.3}
{\em (a)}: $f_{H,\cW} = \frac{\dim V}{|G|} F_{H,\cW}$;\\
{\em (b)}: $F_{H,\cW} = p_H \sum_{i=1}^d a_i \sum_{j=1}^{n_i} g_{ij}$; \\
{\em (c)}: $f_{H,\cW}$ and $F_{H,\cW}$ are elements of the Hecke algebra $\QQ[H \backslash G/H]$.
\end{lem}

\begin{proof}
(a): Using Lemma \ref{lem2.2}, we have,
\begin{eqnarray*}
f_{H,\cW} & = & e_{\cW} p_H \\
& = & \frac{\dim V}{|G||H|} \sum_{i=1}^d \sum_{j=1}^{n_i} \left( \sum_{k \in H} \tr_{K_V/\QQ}(\chi_V ((g_{ij}k)^{-1})) \right)
                           g_{ij}k  \left( \sum_{h \in H} h \right)\\
& = & \frac{\dim V}{|G||H|} \sum_{i=1}^d \sum_{j=1}^{n_i} \left( \sum_{k \in H} \tr_{K_V/\QQ}(\chi_V (kg_{ij}^{-1})) \right)
                           g_{ij}  \left( \sum_{h \in H} h \right)\\
& = & \frac{\dim V}{|G||H|} \sum_{i=1}^d \sum_{j=1}^{n_i} \left( \sum_{k \in H} \tr_{K_V/\QQ}(\chi_V (kg_{i1}^{-1})) \right)
                           g_{ij}  \left( \sum_{h \in H} h \right) = \frac{\dim V}{|G|} F_{H,\cW},\\
\end{eqnarray*}
For the proof of (b) we start with $f_{H,\cW} = p_H e_{\cW}$ and proceed as above. Finally, according to
Lemma \ref{lem2.1} and \eqref{eq2.2} the coefficients of $f_{H,\cW}$ are the same on double cosets, which gives (c).
\end{proof}

\subsection{The integers $a_i$}

The following lemma is a generalization of \cite[91.60 p.391]{McD}, where it is proved under the additional assumption that
$(G,H)$ is a Gelfand pair. Recall that $(G,H)$ is called a Gelfand pair, if $\dim V^H = 0$ or $1$ for every complex irreducible
representation $V$ of $G$. This is equivalent to the fact that the Hecke algebra $\CC[H \backslash G/H]$ is commutative.

Let $V$ be a complex irreducible representation of $G$ and let $( \cdot, \cdot)$ denote any $G$-invariant
hermitian scalar product on $V$ (unique up to a positive constant).

\begin{prop} \label{prop2.4}
Assume $\dim V^H = 1$ for a subgroup $H$ of $G$. Then, for any nonzero vector $v \in V^H$,
$$
(v,gv) = \frac{(v,v)}{|H|} \sum_{h \in H} \chi_V(hg^{-1})
$$
for all $g \in G$.
\end{prop}

\begin{proof}
It follows from Schur's character relations (see e.g. \cite[Proposition 13.6.4]{bl}) that
$$
l_v = \frac{\dim V}{|G|(v,v)} \sum_{g \in G} (v,gv)g
$$
is a primitive idempotent in $\CC[G]$ such that $\CC[G]l_v$ affords the representation $V$.

We claim that $l_v \in \CC[H \backslash G/H]$. In fact, $hl_v = \frac{\dim V}{|G|(v,v)}\sum_{g \in G} (hv,hgv)hg = l_v$, and similarly
$l_vh = l_v$ for all $h \in H$. Therefore $l_v$ is constant on double cosets of $H$ in $G$, which gives the assertion.

Now $l_v \in \CC[G]e_V$ implies $l_v = p_H l_v p_H \in p_H \CC[G]e_Vp_H = f_{H,V} \CC[G] f_{H,V}$. On the other hand,
Lemma \ref{lem2.1} gives
$$
f_{H,V} \CC[G] f_{H,V} \simeq \End_{\CC[G]}(\CC[G]f_{H,V}),
$$
and, as we saw after \eqref{eq2.4}, the left ideal $\CC[G]f_{H,V}$ provides the representation $V$ with multiplicity $\dim V^H = 1$.
Hence Schur's lemma implies that $f_{H,V} \CC[G] f_{H,V}$ is a one-dimensional complex
vector space. Therefore $l_v$ and $f_{H,V}$ differ at most by a constant. So in order to complete the proof,
it suffices to compare the coefficients of $1 \in G$ in $l_v$ and $f_{H,V}$.

Now
$$
f_{H,V} = \frac{\dim V}{|G||H|} \sum _{g\in G} \sum_{h\in H} \chi_V(g^{-1}) gh =
\frac{\dim V}{|G||H|} \sum _{g\in G} \left(\sum_{h\in H} \chi_V(hg^{-1}) \right) g
$$
Hence the coefficient of $1 \in G$ in $f_{H,V}$ is
$$
\frac{\dim V}{|G|} \left( \frac{1}{|H|} \sum_{h \in H} \chi_V(h) \right) = \frac{\dim V}{|G|} \langle V|H,1_H \rangle_H
= \langle V,\rho_H \rangle_G = \frac{\dim V}{|G|} \dim V^H = \frac{\dim V}{|G|}
$$
which coincides with the coefficient of $1 \in G$ of $l_v$. Hence $f_{H,V} = l_v$. Comparing the coefficients in this equation
completes the proof of the proposition.
\end{proof}

As a consequence we obtain the following inequality for the integers $a_i$ of \eqref{eq2.6}.

\begin{prop} \label{prop2.5}
Let the representatives $g_{ij}$ of the cosets of $H$ in $G$ be as above and assume $\dim V^H =1$. Then for
every $1 \leq i \leq d$
the integer $a_i$ satisfies
\begin{equation} \label{eq2.7}
a_i \leq a_1 = [K_V:\QQ] |H|
\end{equation}
with equality if and only if $g_{i1}$ fixes $V^H$.
\end{prop}

\begin{proof}
First consider $a_{1}$. Since we assumed $g_{11} = 1_G$, we have
$$
a_1 = \sum_{h \in H} \tr_{K_V/\QQ}(\chi_V(h)) = \sum_{h \in H} \sum_{\varphi \in \Gal(K_V/\QQ)}\varphi(\chi_V(h)) =
\sum_{\varphi \in \Gal(K_V/\QQ)} \varphi (\sum_{h \in H} \chi_V(h))
$$
But $\sum_{h \in H} \chi_V(h) = |H| \dim V^H = |H|$ and hence $a_1 = [K_V:\QQ] |H|$.

For the proof of the inequality note that the Schur index $m_V$ is a divisor of $\dim V^H = 1$, which implies $L_V = K_V$.
Hence the representation $V$ is defined over $K_V$, i.e. there is a $G$-representation $V_{K_V}$ over $K_V$ such that
$V = V_{K_V} \otimes \CC$.

Consider the $K_V$-vector space
$$
\cV = \bigoplus_{\varphi \in \Gal (K_V/\QQ)} V_{K_V}^{\varphi},
$$
where $V_{K_V}^{\varphi}$ is the representation conjugate to $V_{K_V}$ by $\varphi$. We denote the element of $V_{K_V}^{\varphi}$
corresponding to
$ v \in V_{K_V}$ by $v^{\varphi}$ and observe that
$$
(\alpha v)^{\varphi} = \varphi(\alpha)v^{\varphi}
$$
for all $\alpha \in K_V$. The group $\Gal(K_V/\QQ)$ acts on $\cV$ by permuting the coordinates, and,
if $\Gal(K_V/\QQ) = \{ \varphi_1 = 1, \varphi_2, \ldots, \varphi_{[K_V:\QQ]} \}$,
then its fixed set is the subset of $\cV$ given by
$$
\cW = \{ (v,v^{\varphi_2}, \ldots, v^{\varphi_{[K_V:\QQ]}}) \; | \; v \in V_{K_V} \}.
$$
This is a rational subvector space of $\cV$ of dimension $\dim V \cdot [K_V:\QQ]$, which affords the
rational irreducible representation
defined by $V$. In order to see this, choose a basis $\{v_1, \ldots, v_r \}$ of $V_{K_V}$ and a basis
$\{\eta_1, \ldots, \eta_s \}$ for $K_V/\QQ$. Then $\cW$ is the rational vector space with basis
$$
w_{i,k} = (\eta_kv_i, \ldots, (\eta_kv_i)^{\varphi_s} = \varphi_s(\eta_k) v_i^{\varphi_s}).
$$
Furthermore, it is clear that $\cV \simeq \cW \otimes_{\QQ} K_V$.

Certainly, for any nonzero vector $v \in V_K$ we may choose a $G$-invariant hermitian scalar product $( \cdot, \cdot)$ on $V_{K_V}$
such that $(v,v)$ is a (positive) rational number. Then clearly $S: \cW \times \cW \ra \QQ$, defined by
$$
S(w_1,w_2) = \tr_{K_V/\QQ} (v_1,v_2)
$$
for any $w_i = (v_i, v_i^{\varphi_2}, \ldots, v_i^{\varphi_{[K_V:\QQ]}})$, is a $G$-invariant scalar product on the $\QQ$-vector
space $\cW$. In particular, choosing $(\cdot,\cdot)$ for a fixed nonzero $v \in V_K^H$, and letting
$w = (v, v^{\varphi_2}, \ldots,$  $\varphi^{[K_V:\QQ]}) \in \cW$, we have
by Proposition \ref{prop2.4},
$$
S(w,g_{i1}w) = \tr_{K_V/\QQ}(v,g_{i1}v) = \tr_{K_V/\QQ} \left( \frac{(v,v)}{|H|} \sum_{h \in H} \chi_V(hg_{i1}^{-1}) \right) =
\frac{(v,v)}{|H|} a_i
$$
for $i = 1, \ldots, d$. Since $S$ is positive definite, symmetric and $G$-invariant, we conclude
$$
0 \leq S(g_{i1}w - w, g_{i1}w - w) = 2(S(w,w) - S(w,g_{i1}w)) = \frac{2(v,v)}{|H|}( a_1 - a_i)
$$
with equality if and only if $g_{i1}w = w$, which completes the proof of the proposition.
\end{proof}

\section{The correspondences}\label{S:3}

\subsection{The set up}
Let $Z$ be a smooth projective curve such that the group $G$ acts on $Z$. Let $JZ$ denote the Jacobian of $Z$ and
$\End_{\QQ}(JZ) = \End(JZ) \otimes_{\ZZ} \QQ$ its endomorphism algebra. The group action induces an algebra homomorphism
$$
\QQ[G] \ra \End_{\QQ}(JZ).
$$
Since this homomorphism is canonical, we will denote the elements of $\QQ[G]$ and their images by the same letter.
In particular we consider elements of $\ZZ[G]$ as endomorphisms of $JZ$. For any $\alpha \in \QQ[G]$ we define
$$
\im(\alpha) := \im (m \alpha) \subset JZ,
$$
where $m$ is a positive integer such that $m \alpha \in \End(JZ)$. It is an abelian subvariety of $JZ$ which
certainly does not depend on the chosen integer $m$.

Now consider a subgroup $H$ of $G$. If we denote the quotients of $Z$ by $H$ and $G$ by $X = Z/H$ and $Y = Z/G$ respectively,
we have the following diagram
\begin{equation} \label{diag} \xymatrix{ Z \ar[dd]_{\pi}
\ar[dr]^{\pi_H} \\ & X \ar[dl]^{\varphi_H} \\ Y }
\end{equation}

Since $\QQ[H \backslash G/H] = p_H \QQ[G]  p_H$, the homomorphism $\QQ[G] \ra \End_{\QQ}(JZ)$ restricts
to an algebra homomorphism
$$
\QQ[H \backslash G/H] \ra \End_{\QQ}(p_H(JZ)).
$$
Now the pull-back map $\pi_H^*: JX \ra p_H(JZ)$ and the restriction of the norm map $\Nm_{\pi_H}: p_H(JZ) \ra JX$  are isogenies
satisfying $\Nm \pi_H \circ \pi_H^* = |H| 1_{JX}$. This implies that the composition
\begin{equation} \label{eq3.2}
\varepsilon: \QQ[H \backslash G/H] \ra \End_{\QQ}(JX), \quad \varphi \mapsto \frac{1}{|H|} \Nm_{\pi_H} \circ \varphi \circ \pi_H^*
\end{equation}
is a homomorphism of $\QQ$-algebras.\\

\subsection{The Hecke ring $\frac{1}{|H|} \ZZ[H \backslash G/H]$}
Let $H_1 = H, H_2, \ldots, H_d$ denote the double cosets of
$H$ in $G$. Consider for $i = 1, \ldots, d$ the elements of $\QQ[H \backslash G/H]$ defined by
$$
F_i := \frac{1}{|H|} \sum_{g \in H_i} g = \sum_{j=1}^{n_i} g_{ij} p_H = p_H \sum_{j=1}^{n_i} g_{ij}
$$
with $p_H$ as in \eqref{eq2.3} and representatives $g_{ij}$ as chosen in subsection 2.2.
We define
$$
\frac{1}{|H|} \ZZ[H \backslash G/H]
$$
as the free $\ZZ$-module with basis $F_1, \ldots, F_d$.
So any element
$F \in \frac{1}{|H|} \ZZ[H\backslash G/H]$ can be written as
$$
F = \sum_{i=1}^d \alpha_i F_i
= \sum_{i=1}^d \sum_{j=1}^{d_i} \alpha_i g_{ij} p_H = p_H \sum_{i=1}^d \alpha_i \sum_{j=1}^{d_i} g_{ij}
$$
with uniquely determined integer constants $\alpha_i$. The following lemma justifies the notation.
\begin{lem} \label{lem3.1}
$\frac{1}{|H|} \ZZ[H \backslash G/H]$ is a $\ZZ$-algebra of rank $d$ with unit element $F_1 = p_H$.
\end{lem}

\begin{proof}
It remains to show that $F_iF_j \in \frac{1}{|H|} \ZZ[H \backslash G/H]$ for all $i$ and $j$. Now it is easy to see that
$$
F_iF_j = \sum_{k=1}^d c_{ijk} F_k \quad \mbox{with} \quad  c_{ijk} = \frac{|H_i \cap g_kH_j^{-1}|}{|H|}
$$
where $g_k$ is any element of $H_k$. Since both $H_i$ and $g_kH_j^{-1}$
are unions of left cosets of $H$ in $G$, the constants $c_{ijk}$ are integers which implies the assertion.
\end{proof}
We call $\frac{1}{|H|} \ZZ[H \backslash G/H]$ {\it the Hecke ring of the subgroup} $H$ of $G$.

\subsection{The homomorphism $\frac{1}{|H|} \ZZ[H \backslash G/H] \ra \Div(X \times X)$} It is well known that
the canonical homomorphism $\QQ[G] \ra \End_{\QQ}(JZ)$ of Section 3.1 factorizes via the algebra of $\QQ$-correspondences $\Div_{\QQ}(Z \times Z)$.
In this subsection we show that the homomorphism $\epsilon$ of \eqref{eq3.2} factorizes via the algebra $\Div_{\QQ}(X \times X)$.
To see this, note first that every element $F = \sum_{i=1}^d \alpha_i F_i \in \QQ[H \backslash G/H]$ induces a $\QQ$-correspondence on $Z$ defined by
$$
\cD_F = \frac{1}{|H|} \sum_{i=1}^d \alpha_i \sum_{j=1}^{n_i} \sum_{h \in H} \Gamma_{g_{ij}h}
= \frac{1}{|H|} \sum_{i=1}^d \alpha_i \sum_{j=1}^{n_i} \Gamma_{g_{ij}} \sum_{h \in H} \Gamma_{h}  \in \Div_{\QQ}(Z \times Z)
$$
where $\Gamma_g \subset Z \times Z$ denotes the graph of $g$ for any $g \in G$, and we use the multiplication of the
ring $\Div_{\QQ}(Z \times Z)$. The corresponding map
$Z \ra \Div_{\QQ} Z$ is given by
$$
\cD_F(z) = \frac{1}{|H|} \sum_{i=1}^d \alpha_i \sum_{j=1}^{n_i} \sum_{h \in H} g_{ij} h (z) =
\frac{1}{|H|} \sum_{i=1}^d \alpha_i \sum_{h \in H} \sum_{j=1}^{n_i} h g_{ij}(z),
$$
since the $g_{ij}$ are representatives for both the left and right cosets of $H$ in $G$.
In particular we have for every $h \in H$ and $z \in Z$,
\begin{equation} \label{eq3.3}
\cD_F(h(z)) = \cD_F(z).
\end{equation}
As every $\QQ$-correspondence on $Z$,
the correspondence $\cD_F$ pushes down to a correspondence on $X$, namely
$$
{\overline \cD}_F := \frac{1}{|H|} (\pi_H \times \pi_H)_* \cD_F \in \Div_{\QQ}(X \times X).
$$ The
following proposition shows that
${\overline \cD}_F$ is actually an integral correspondence for any $F \in \frac{1}{|H|} \ZZ[H \backslash G/H]$.

\begin{prop} \label{prop3.1}
The map $F \mapsto {\overline \cD}_F$ is a homomorphism of rings
$$
\frac{1}{|H|} \ZZ[H\backslash G/H] \ra \Div(X \times X).
$$
If $F = \sum_{i=1}^d \alpha_i F_i$, we have for any $x \in X$
\begin{equation} \label{eq3.4}
{\overline \cD}_F(x) = \sum_{i=1}^d \alpha_i \sum_{j=1}^{n_i} \pi_H(g_{ij}(z))
\end{equation}
where $z \in Z$ is a preimage of $x$.
\end{prop}

\begin{proof}
First note that the right hand side of \eqref{eq3.4} is independent of the choice of the preimage $z$ of $x$.
If $z'$ is another preimage, there is an $h \in H$ such that $z' = h(z)$. The choice of the $g_{ij}$ as a simultaneous set of
representatives for the left and right cosets of $H$ in $G$ implies for any $i, \; 1 \leq i \leq d$,
$$
\sum_{j=1}^{n_i} \pi_H (g_{ij}(z')) = \sum_{j=1}^{n_i} \pi_H (g_{ij}(h(z))) = \sum_{j=1}^{n_i} \pi_H (h(g_{ij}(z)))
= \sum_{j=1}^{n_i} \pi_H (g_{ij}(z)).
$$
Multiplying by $\alpha_i$ and summing over all $i$ gives the assertion.

Hence the right hand side of \eqref{eq3.4} defines a set-theoretical correspondence ${\overline \cD}_F$. In order to show that
${\overline \cD}_F$ is an algebraic correspondence on $X$ and equals $\frac{1}{|H|} (\pi_H \times \pi_H)_*\cD_F$, it suffices to show that
$((\pi_H \times \pi_H)_* \cD_F) (x)= |H| {\overline \cD}_F (x)$ for all $x \in X$.

Using the special properties of the representatives $g_{ij}$ and choosing a fixed preimage $z$ of $x$, we have
\begin{eqnarray*}
((\pi_H \times \pi_H)_*\cD_F) (x) = \sum_{h \in H} \pi_H(\cD_F(h(z))) & = & |H| \cdot \pi_H(\cD_F(z))\\
& = & \sum_{h \in H} \sum_{i=1}^d \alpha_i \sum_{j=1}^{n_i} \pi_H(hg_{ij}(z))\\
& = & |H| \cdot \sum_{i=1}^d \alpha_i \sum_{j=1}^{n_i} \pi_H(g_{ij}(z)) = |H| \cdot {\overline \cD}_F(x).
\end{eqnarray*}
It remains to see that the map $F \mapsto {\overline \cD}_F$ is a ring-homomorphism. But it is easily seen using
\eqref{eq3.3} several times that ${\overline \cD}_{FF'}(x) = {\overline \cD}_{F} {\overline \cD}_{F'}(x)$ for any
$F, F' \in \frac{1}{|H|} \ZZ[H\backslash G/H]$ and any $x \in X$.
\end{proof}

Since $F_1, \ldots, F_d$ are also a $\QQ$-basis of the Hecke algebra $\QQ[H \backslash G/H]$, the homomorphism of Proposition
\ref{prop3.1} extends to a homomorphism of $\QQ$-algebras $\QQ[H \backslash G/H] \ra \Div_{\QQ}(X \times X)$.
Together with the definition
of the homomorphism $\epsilon: \QQ[H \backslash G/H] \ra \End_{\QQ}(JX)$, we obtain

\begin{cor} \label{cor3.3}
The homomorphism of Proposition \ref{prop3.1} extends to a homomorphism of $\QQ$-algebras
$\QQ[H \backslash G/H] \ra \Div_{\QQ}(X \times X)$,
which factorizes the homomorphism  $\epsilon: \QQ[H \backslash G/H] \ra \End_{\QQ}(JX)$.
\end{cor}

Consider now a rational irreducible representation $\cW$ of $G$ with associated complex irreducible representation $V$.
From subsection \ref{ss2.2} we conclude that
$$
F_{H,\cW} = \sum_{i=1}^d a_i F_i \in \frac{1}{|H|} \ZZ[H \backslash G/H]
$$
with integers $a_i$ given by \eqref{eq2.6}. Hence we may apply Proposition \ref{prop3.1} to $F_{H,\cW}$ to conclude
the first part of the following proposition.

\begin{prop} \label{prop3.2}
${\overline \cD}_{H,\cW} := {\overline \cD}_{F_{H,\cW}}$ is an integral symmetric correspondence on $X$, of degree 0 if $\cW$
is non-trivial and of
degree $|G|$ if $\cW$ is trivial.
\end{prop}

\begin{proof}
According to Proposition \ref{prop3.1}, ${\overline \cD}_{H,\cW}$ is an integral correspondence on $X$.
Applying \eqref{eq2.6}, \eqref{eq2.2} and Lemma \ref{lem2.2} we have for the degree,
$$
\begin{array}{lll}
\deg {\overline \cD}_{H,\cW} &
= {\displaystyle \sum_{i=1}^d a_in_i  = \sum_{i=1}^d n_i \sum_{h\in H} \tr_{K_V/\QQ}(\chi_V(hg_{i1}^{-1}))} \\ \\
& = {\displaystyle \sum_{g \in G} \tr_{K_V/\QQ}(\chi_V(g))} \\ \\
& =\left\{ \begin{array}{cll}
                                                 0 & \; \mbox{if} & \cW \; \mbox{is non-trivial;}\\
                                                 |G| & \; \mbox{if} & \cW \; \mbox{is trivial}.
                                                 \end{array} \right.
\end{array}
$$
where the zero in the last equation is just the fact that the representation $\cW$ and the trivial representation of $G$
have zero scalar product.

It remains to show that the correspondence ${\overline \cD}_{H,\cW}$ is symmetric. For any $x \in X$ choose a preimage $z \in Z$.
Then $\pi_H(g_{ij}(z))$ appears with multiplicity $a_i$ in ${\overline \cD}_{H,\cW}(x)$.
It suffices to show that $x$ appears with the
same multiplicity in ${\overline \cD}_{H,\cW}(\pi_H(g_{ij}(z)))$.

By definition we have
$$
{\overline \cD}_{H,\cW}(\pi_H(g_{ij}(z))) =  \sum_{k=1}^d a_k \sum_{l=1}^{n_k} \pi_H(g_{kl}(g_{ij}(z))).
$$
We are interested in its coefficient of $x = \pi_H(z)$. But $\pi_H(g_{kl}(g_{ij}(z))) = x$ if and only if $g_{kl}g_{ij} \in H$,
which is the case if and only if
\begin{equation*}
g_{kl} \in Hg_{ij}^{-1}
\end{equation*}
Hence $Hg_{kl}H = Hg_{ij}^{-1}H$
and it suffices to show that for any $g \in G$,
$F_{H,\cW}$ has the same coefficient on the double cosets $HgH$ and $Hg^{-1}H$.
But this follows from equation \eqref{eq2.2} and
\eqref{eq2.6} or, to be more precise, from the fact that $\tr_{K_V/\QQ}(\chi_V(g)) = \tr_{K_V/\QQ}(\chi_V(g^{-1}))$.
This completes the proof
of the proposition.
\end{proof}

\begin{prop} \label{prop3.3} For any two different rational irreducible representations $\cW$ and ${\widetilde \cW}$ of $G$ we have

{\em (a):} \hspace{5cm} ${\overline \cD}_{H,\cW}^2 = \frac{|G|}{\dim V} {\overline \cD}_{H,\cW},$\\

{\em (b):} \hspace{5cm} ${\overline \cD}_{H,\cW} {\overline \cD}_{H,{\widetilde \cW}} = 0$.
\end{prop}

\begin{proof}
(a) follows from Lemma \ref{lem2.3} and Proposition \ref{prop3.1} together with the fact that
$f_{H.\cW} = \frac{|G|}{\dim V}F_{H,\cW}$ is an idempotent. Similarly, (b) is a consequence of the same results together with
the fact that
$f_{H,\cW}$ and $f_{H,{\widetilde \cW}}$ are orthogonal idempotents.
\end{proof}

\subsection{The trace correspondence} We denote by $\cT_{X/Y} := \varphi_H^* \varphi_H$
the trace correspondence of the
covering $\varphi_H: X \ra Y$. So for any $x \in X$,
$$
\cT_{X/Y}(x) = \sum_{i=1}^d \sum_{j=1}^{n_i} \pi_H(g_{ij}(z)).
$$
where $z \in Z$ is a preimage of $x$. On the other hand, if $V_0 = \cW_0$ denotes the trivial representation of $G$,
according to \eqref{eq3.4} and the proof
of Proposition \ref{prop3.2} we have
$$
{\overline \cD}_{H,V_0}(x) = \sum_{i=1}^d |H| \sum_{j=1}^{n_i} \pi_H(g_{ij}(z))
$$
which implies
$$
\cT_{X/Y} = \frac{1}{|H|} {\overline \cD}_{H,V_0}.
$$
Hence as a special case of Proposition \ref{prop3.3} we obtain
\begin{cor} \label{lem3.5}
{\em (a)} \hspace{2.6cm} $\cT_{X/Y}^2 = [G:H] \cT_{X/Y}$;

{\em (b)} For any nontrivial rational irreducible representation $\cW$ of $G$,
$$
{\overline \cD}_{H,\cW} \cT_{X/Y} = \cT_{X/Y} {\overline \cD}_{H,\cW} = 0.
$$
\end{cor}

\subsection{The Kanev correspondence}
Let $\cW_1, \ldots, \cW_r$ denote nontrivial pairwise non-isomorphic rational irreducible representations of the group $G$
with associated complex irreducible representations $V_1, \ldots, V_r$.
We make the following hypothesis for the $\cW_i$ and a subgroup $H$ of $G$.

\begin{hyp} \label{hyp}
For all $k,l = 1 ,\ldots, r$ we assume\\
{\em a)} \quad $\dim V_k = \dim V_l =: n$ fixed, \\
{\em b)} \quad $K_{V_k} = K_{V_l} =: L$ fixed, \\
{\em c)} \quad $\dim V_k^H = 1$,\\
{\em d)} \quad $H$ is maximal with property c), that is, for every subgroup $N$ of $G$ with $H \subsetneqq N$
there \hspace*{0.8cm} is an index $k$ such that $\dim V_k^N = 0$.
\end{hyp}

Recall the correspondence ${\overline \cD}_{H,\cW_k}$ given by
\begin{equation} \label{e3.5}
{\overline \cD}_{H,\cW_k}(x) = \sum_{i=1}^d a_{ki} \sum_{j=1}^{n_i} \pi_H g_{ij}(z)
\end{equation}
for all $x \in X$ and $z \in Z$ with $\pi_H(z) = x$, and where the $a_{ki}$ are the integers
$$
a_{ki} = \sum_{h \in H} \tr_{L/\QQ}(\chi_{V_k}(hg_{i1}^{-1})).
$$
Now consider the correspondence
$$
{\overline \cD} := \sum_{k=1}^r {\overline \cD}_{H,\cW_k}.
$$
Denoting $b_i := \sum_{k=1}^r a_{ki}$ for $i = 1, \ldots, d$ we have
$$
{\overline \cD}(x) = \sum_{i=1}^d b_{i} \sum_{j=1}^{n_i} \pi_H g_{ij}(z).
$$
Then set
\begin{equation} \label{eqn3.5}
b := \gcd \{ b_1 - b_i \; | \; 2 \leq i \leq d \}.
\end{equation}
Denoting by $\Delta_X \in \Div(X \times X)$ the identity correspondence, we define the Kanev correspondence on $X$
associated to $\cW_1, \ldots , \cW_r$ (which we omit in the notation) by
\begin{equation} \label{eq3.5}
\cK_X := \Delta_X - \frac{1}{b} {\overline \cD} + (\frac{b_1}{b} - 1 ) \cT_{X/Y}.
\end{equation}

\begin{lem} \label{lem3.6}
Suppose $\cW_1, \ldots, \cW_r$ are nontrivial pairwise non-isomorphic rational irreducible representations of $G$
satisfying the Hypothesis \ref{hyp} and recall $L= K_{V_i}$ for all $i$.
Then $\cK_X$ is an integral effective symmetric correspondence on $X$ of degree
\begin{equation} \label{eq3.6}
\deg \cK_X = 1 + \frac{r|G|}{b} [L:\QQ] - [G : H].
\end{equation}
\end{lem}

\begin{proof}
Since for all $k$, $V_k$ is not the trivial representation, but satisfies $\dim V_k^H = 1$, it follows from Proposition
\ref{prop2.5} that $a_{ki} \leq a_{k1}$ for all $i,k$. Therefore $b_i \leq b_1$ for all $i$.
Moreover, the maximality condition d) of $H$ together with Proposition \ref{prop2.5} imply that actually
$b_i < b_1$ for all $i \geq 2$. Therefore $b$ is a positive integer.

Now choose for every $x \in X$ a preimage $z \in Z$. Then
\begin{eqnarray*}
\cK_X(x) & = & x - \frac{1}{b} \sum_{i=1}^d b_{i} \sum_{j=1}^{n_i} \pi_H (g_{ij}(z)) +
(\frac{b_1}{b} -1) \sum_{i=1}^d \sum_{j=1}^{n_i} \pi_H(g_{ij}(z))\\
& = & x + \sum_{i=1}^d \left(\frac{b_1 - b_i}{b} -1\right) \sum_{j=1}^{n_i} \pi_H(g_{ij}(z))
= \sum_{i=2}^d \left(\frac{b_1 - b_i}{b} -1\right) \sum_{j=1}^{n_i} \pi_H(g_{ij}(z)).
\end{eqnarray*}
Hence $\cK_X$ is integral and effective. The symmetry of $\cK_X$ is a consequence of Proposition \ref{prop3.2}. Finally,
using that $\deg {\overline \cD} = 0$ according to Proposition \ref{prop3.2} as well as Proposition \ref{prop2.5}, we have
%using \eqref{eq2.7}, \eqref{eq2.6} and \eqref{eq2.2},
\begin{eqnarray*}
\deg \cK_X & = & 1 + (\frac{b_1}{b} -1) \deg \cT_{X/Y}\\
& = & 1 + (\frac{r[L:\QQ]|H|}{b} - 1) [G:H] = 1 + \frac{r|G|}{b}[L: \QQ] - [G:H].
\end{eqnarray*}
%\begin{eqnarray*}
%\deg \cK_X & = & 1 + \sum_{i=1}^d (\frac{b_1 - b_i}{b} - 1)n_i\\
%& = & 1 + \frac{r|H|}{b}[K_1: \QQ] \sum_{i=1}^d n_i -\frac{1}{b}\sum_{i=1}^d \sum_{k=1}^r a_{ki}n_i - \sum_{i=1}^d n_i\\
%& = & 1 + \frac{r|H|}{b}[K_1 : \QQ][G : H] - \frac{1}{b} \sum_{k=1}^r \tr_{K_{V_k}/\QQ}(\sum_{g \in G} \chi_V(g)) - [G:H]\\
%& = & 1 + \frac{r|G|}{b}[K_1:\QQ] - [G:H]
%\end{eqnarray*}
%where the last equation holds, because $\sum_{g \in G} \chi_{V_k}(g) = 0$, since $V_k$ is not the trivial representation for all
%$k = 1, \ldots r$.
\end{proof}

\begin{prop} \label{prop3.9}
Suppose $\cW_1, \ldots, \cW_r$ are nontrivial pairwise non-isomorphic rational irreducible representations of $G$
satisfying the Hypothesis \ref{hyp}.
Then the correspondence $\cK_X$ satisfies the cubic equation
\begin{equation} \label{eq3.7}
(\cK_X - \Delta_X)(\cK_X + (q-1)\Delta_X)(\cK_X - \deg \cK_X \cdot \Delta_X) = 0
\end{equation}
where $\Delta_X$ is the diagonal in $X \times X$ and $q$ the positive integer
\begin{equation} \label{eq3.8}
q = \frac{|G|}{b \cdot n}.
\end{equation}
\end{prop}

\begin{proof}
Note first that Corollary \ref{lem3.5} and \eqref{eq3.6} imply
\begin{equation} \label{eq3.9}
\cK_X \cT_{X/Y} = \left(1 + (\frac{b_1}{b} -1)[G:H]\right) \cT_{X/Y} = \deg \cK_X \cdot \cT_{X/Y}.
\end{equation}
Now using the definition of $\cK_X$, Proposition \ref{prop3.3}, Corollary \ref{lem3.5} and equation \ref{lem3.5} we get
\begin{eqnarray*}
(\Delta_X - \cK_X)^2 & = & (\frac{1}{b} {\overline \cD} - (\frac{b_1}{b} -1) \cT_{X/Y})^2\\
& = & \frac{1}{b}{\overline \cD}^2 + (\frac{b_1}{b} - 1)\cT_{X/Y}^2\\
& = & \frac{|G|}{b^2 \cdot n}{\overline \cD} + (\frac{b_1}{b} - 1)^2 [G:H] \cT_{X/Y}\\
\end{eqnarray*}
Defining $q$ by the right hand side of \eqref{eq3.8}, this gives
\begin{equation}  \label{eq3.10}
(\cK_X - \Delta_X)^2 + q (\cK_X - \Delta_X) + c \cdot \cT_{X/Y} = 0
\end{equation}
where $c$ denotes the rational number
$$
c = (1 - \frac{b_1}{b})\left((\frac{b_1}{b} -1)[G:H] + \frac{|G|}{b \cdot n} \right).
$$
Multiplying \eqref{eq3.10} by $\cK_X - \deg \cK_X \cdot \Delta_X$ and applying \eqref{eq3.9}, we get
\eqref{eq3.7}.

It remains to show that $q$ is an integer. For this consider a general point $y \in Y$. The action of $\cK_X$,
respectively  of $\cT_{X/Y}$, on the fibre $\varphi_H^{-1}(y)$ is described by a square  integral matrix of size
$[G:H]$ denoted by
$M_{\cK_X}$, respectively by $M_{\cT_{X/Y}}$. If we denote $N = M_{\cK_X} - E$, where $E$ denotes the identity
matrix of size $[G:H]$, equation \eqref{eq3.10}
implies
\begin{equation} \label{eq3.11}
N(N +qE) = -cM_{\cT_{X/Y}}.
\end{equation}
We will complete the proof by showing that the rational number $-q$ is an
eigenvalue of the integral matrix $N$.
If $c =0$, then equation \eqref{eq3.11} implies that $N$ satisfies the
equation $x(x+q)=0$.
Since $N$ is not the zero matrix, then either $x+q$ or $x(x+q)$ is the
minimal polynomial for $N$ and we are done.
If $c \neq 0$, then, according to \eqref{eq3.7} and \eqref{eq3.11}, the
minimal polynomial of $N$ is $x(x + q)(x-(\deg \cK_X -1))$
and we conclude in the same way that $q$ is an integer.\end{proof}

\begin{rem}
Equation \eqref{eq3.11} implies that $c M_{\cT_{X/Y}}$ is an integral matrix. Since all entries of the
matrix $M_{\cT_{X/Y}}$ are equal to 1, we conclude that $c$ is also an integer.
\end{rem}

\begin{cor}
With the notation of Proposition \ref{prop3.9} the following conditions are equivalent:\\
{\em (a):} $q=1$,\\
{\em (b):} $\rho_H \simeq \cW_0 \oplus \cW_1 \oplus \cdots \oplus \cW_r$,\\
{\em (c):} $\cK_X = 0$.
\end{cor}

\begin{proof}
Note first that
$$
\chi_{\rho_H}(1_G) = [G:H] \geq 1 + r n [L:\QQ].
$$
According to Proposition \ref{prop3.9} $q=1$ if and only if $|G| = b n$. Since $\cK_X$ is effective, this is equivalent to
$$
0 \leq \deg \cK_X = 1 + r n [L:\QQ] - [G:H] \leq 0.
$$
This implies all assertions.
\end{proof}

\section{Prym varieties}\label{S:4}

\subsection{The Prym variety $P_{\overline \cD}$}\label{SS:prym}
Let the notation be as in the last section. So $H$ is a subgroup of $G$. Moreover denote by $\cW := (\cW_1, \ldots, \cW_r)$
an $r$-tuple of nontrivial pairwise non-isomorphic rational irreducible representations of the group $G$
with associated complex irreducible representations $V_1, \ldots, V_r$ satisfying Hypothesis \ref{hyp}.
Let $\delta_{\overline \cD}$ denote the endomorphism of the Jacobian $JX$ associated to the correspondence ${\overline \cD}$.
Then we get the following statement, as a direct consequence of Proposition  \ref{prop3.3} and the definitions.

\begin{prop} \label{prop4.1}
$$
\delta_{\overline \cD}^2 = \frac{|G|}{n} \delta_{\overline \cD}.
$$
\end{prop}

We denote by
$$
P_{\overline \cD} := \im(\delta_{\overline \cD})
$$
the image of the endomorphism $\delta_{\overline \cD}$ in the Jacobian $JX$ and call it the {\it (generalized) Prym variety}
associated to the correspondence ${\overline \cD}$. Our aim is to investigate the restriction of the canonical polarization
of $JX$ to $P_{\overline \cD}$ in some cases.\\

Let us denote by $\tau_{X/Y} \in \End(JX)$ the endomorphism associated to the trace correspondence $\cT_{X/Y}$. The {\it Prym
variety of the covering} $\varphi_H: X \ra Y$ is defined by
$$
\Prym(X/Y) := \ker(\tau_{X/Y})^0 \subset JX.
$$
According to Corollary \ref{lem3.5} we have $\tau_{X/Y}\delta_{\overline \cD} = 0$. This implies

\begin{prop} \label{prop4.3}
The Prym variety $P_{\overline \cD}$ is a subvariety of the Prym variety of the covering $\varphi_H$:
$$
P_{\overline \cD}  \subset \Prym(X/Y).
$$
\end{prop}

Let $\kappa_X$ denote the endomorphism of $JX$ associated to the correspondence $\cK_X$. According to
Corollary \ref{lem3.5} we have
$\kappa_X \tau_{X/Y} = \tau_{X/Y} \kappa_X =(1 + (\frac{b_1}{b}-1)[G:H])\tau_{X/Y}$. Hence $\kappa_X$ restricts to an endomorphism ${\tilde \kappa}_X$ of $\Prym(X/Y)$.

\begin{prop} \label{prop4.4}
{\em (a)}  \hspace{2cm} $P_{\overline \cD} = \im (1_{\Prym(X/Y)} - {\tilde \kappa}_X)$;\\
{\em (b)}  ${\tilde \kappa}_X$ satisfies the quadratic equation
\begin{equation} \label{eq4.1}
{\tilde \kappa}_X^2 + (q - 2){\tilde \kappa}_X - (q-1)1_{\Prym(X/Y)} = 0
\end{equation}
with $q = \frac{|G|}{b \cdot n}$.
\end{prop}

\begin{proof}
According to \eqref{eq3.5} we have $\delta_{\overline D}|_{\Prym(X/Y)} = b( 1_{\Prym(X/Y)} - {\tilde \kappa}_X)$, which gives (a). Finally, (b) is a consequence of \eqref{eq3.10}.
\end{proof}

\begin{prop} \label{prop4.5}
\begin{equation} \label{eq4.2}
\dim P_{\overline \cD} = \frac{1}{q}\left( g(X) + [G:H] - 1 - \frac{r|G|}{b}[L:\QQ] + \frac{1}{2}(\cK_X \cdot \Delta_X ) \right).
\end{equation}
\end{prop}

\begin{proof}
According to Proposition \ref{prop4.1} the element $\frac{n}{|G|}\delta_{\overline \cD}$
is the symmetric idempotent corresponding to the abelian
subvariety $P_{\overline \cD}$ of $JX$. Hence \cite[Corollary 5.3.10]{bl} gives
$$
\dim P_{\overline \cD} = \tr_a \left( \frac{n}{|G|} \delta_{\overline \cD} \right) = \frac{b \cdot n}{|G|} \tr_a(1_{JX} - \kappa_X) =
 \frac{1}{q} \left( g(X) - \frac{1}{2} \tr_r(\kappa_X) \right).
$$
On the other hand, according to \cite[Proposition 11.5.2]{bl} we have for the rational trace of $\kappa_X$,
$$
\tr_r(\kappa_X) = 2 \deg \cK_X - (\cK_X \cdot \Delta_X).
$$
So \eqref{eq3.6} gives the assertion.
\end{proof}

\subsection{Variation of $H$}
Up to now we considered a fixed subgroup $H$ of $G$. For any
$r$-tuple of nontrivial pairwise non-isomorphic rational irreducible representations
$\cW := (\cW_1, \ldots, \cW_r)$ of the group $G$
with associated complex irreducible representations $V_1, \ldots, V_r$ satisfying Hypothesis \ref{hyp}, we associated a
correspondence ${\overline \cD}$ of the curve $X$ and an abelian subvariety of $JX$, the Prym variety $P_{\overline \cD}$.

Suppose now, we are given two subgroups $H_1$ and $H_2$ of $G$ such that $\cW$ satisfies Hypothesis \ref{hyp} for both of them.
For $i = 1, 2$ denote by $X_i$ the curve $Z/H_i$, and by ${\overline \cD}_i$ the associated correspondence.
The following proposition shows how the corresponding Prym varieties $P_{{\overline \cD}_i}$ are related.

\begin{prop}\label{prop:conjugated}
{\em (a):} With the above notations the abelian varieties $P_{{\overline \cD}_1}$ and $P_{{\overline \cD}_2}$ are isogenous.\\
{\em (b):} If in addition $H_1$ and $H_2$ are conjugate in $G$, the canonical isomorphism $X_1 \ra X_2$ induces an
isomorphism of polarized abelian varieties $P_{{\overline \cD}_1} \ra P_{{\overline \cD}_2}$
\end{prop}

\begin{proof}
We only give a sketch, since we do not need the result in the sequel. First note that Proposition \ref{prop3.3} (b) reduces the proof to
the case of one representation, i.e. to $\cW = \cW_1$. For (a) let $f$ denote a primitive idempotent of $\QQ[G]e_{\cW}$
considered as an element of $\End_{\QQ}(JZ)$. Then it is easy to see that $P_{{\overline \cD}_i}$ is isogenous to
$\im(f) \subset JZ$. The proof of (b) is a straightforward computation.
\end{proof}

\begin{rem}
An example of Proposition \ref{prop:conjugated} with $H_1$ and $H_2$ not conjugate but still each one satisfying Hypothesis \ref{hyp}, is provided by the Alternating group of degree $4$, with $V$ the standard representation and the subgroups $H_1$ of order $2$ and $H_2$ of order $3$.
\end{rem}

\subsection{The case $r=1$}
Let $\cW$ denote a rational irreducible representation of $G$.
Recall from \cite{caro} that for any rational irreducible representation $\cW$ and any subgroup $H$ of $G$ there is a uniquely
determined abelian subvariety of the Jacobian $JX$, called the {\it isotypical component} associated to $\cW$ (see also \cite{lre}).
It is isogenous to $B_{\cW}^{\frac{\dim V^H}{m_V}}$, where $B_{\cW}$ is the image of a primitive idempotent of the group algebra $\QQ[G]$ in $JX$
corresponding to a minimal left ideal of the simple subalgebra of $\QQ[G]$ defined by $\cW$. The abelian subvariety $B_{\cW}$
is only
determined up to isogeny in general.

\begin{prop} \label{prop4.2}
Let $\cW$ denote the rational irreducible representation of $G$ associated to $V$ with $\dim V^H = 1$. Then
$P_{{\overline \cD}_{\cW}} = B_{\cW}$ is the isotypical component of $JX$ associated to $\cW$.
\end{prop}

\begin{proof}
This is a consequence of \cite[Proposition 5.2]{caro}: On the one hand $P_{{\overline \cD}_{\cW}}$ is the isotypical component
associated to $\cW$ in $JX$ and on the other hand
it is shown that $P_{{\overline \cD}_{\cW}}$ is isogenous to $B_{\cW}^{\frac{\dim V^H}{m_V}}$.
The assumption $\dim V^H =1$ implies $m_V = 1$ and thus the assertion.
\end{proof}

\subsection{Fixed points of $\cK_X$}
As before, $\cW_1, \ldots, \cW_r$
denote nontrivial pairwise non-isomor\-phic rational irreducible representations of the group $G$
with associated complex irreducible representations $V_1, \ldots, V_r$ satisfying Hypothesis \ref{hyp}.
The number of fixed points of the correspondence $\cK_X$ is by
definition the intersection number $(\cK_X \cdot \Delta_X)$. It depends
on the type of the Galois covering $\pi: Z \ra Y$. In order to express it in terms of the type of $\pi$,
we use the notion of geometric signature as defined
in the introduction.

\begin{prop} \label{prop4.6}
Suppose the action of the group $G$ on the curve $Z$ has geometric signature $[0;(C_1,m_1), \ldots , (C_t,m_t)]$. Then the number of fixed points
of the correspondence $\cK_X$ of the curve $JX$ is given by
$$
(\cK_X \cdot \Delta_X) =  \sum_{j=1}^t m_j \left( q [L:\QQ]\sum_{i=1}^r (\dim V_i - \dim V_i^{G_j})  -
([G:H] -|H \backslash G/G_j|) \right).
$$
where $G_j$ denotes any subgroup in the class $C_j$.
\end{prop}

\begin{proof}
According to \cite[Corollary 3.4]{R} we have for the genus of $X$,
$$
g(X) = 1 - [G:H] + \frac{1}{2} \sum_{j=1}^t m_j([G:H] - |H \backslash G/G_j|).
$$
Since the representations $\cW_i$ are pairwise non-isomorphic, we certainly have
$$
\dim P_{\overline \cD} = \sum_{i=1}^r \dim P_{\cD_{\cW_i}}.
$$
where $P_{\cD_{\cW_i}} = \im(\delta_{\cW_i})$ and $\delta_{\cW_i}$ is the element of $\End_{\QQ}(JX)$ associated
to the correspondence $\cD_{{\overline \cW}_i}$. Now \cite[Theorem 5.12]{R} gives, using Proposition \ref{prop4.2} and the
fact that the Schur index of $V_i$ is 1,
$$
\dim P_{\cD_{\cW_i}} = [L:\QQ] \left( \frac{1}{2} \sum_{j=1}^t m_j(\dim V_i - \dim V_i^{G_j}) -  n \right).
$$
Hence
$$
\dim P_{\overline \cD} = [L:\QQ] \sum_{i=1}^r \left( \frac{1}{2}\sum_{j=1}^t m_j(\dim V_i -
\dim V_i^{G_j}) -  n \right).
$$
Inserting these equations into \eqref{eq4.2} we get
\begin{eqnarray*}
(\cK_X \cdot \Delta_X) & = & 2 \left( \frac{|G|}{b \cdot n} \dim P_{\overline \cD} - g(X) - [G:H] + 1 + \frac{r|G|}{b} [L:\QQ] \right)\\
& = & \sum_{j=1}^t m_j \left( q[L:\QQ]\sum_{i=1}^r  (\dim V_i - \dim V_i^{G_j})  -
([G:H] -|H \backslash G/G_j|) \right).
\end{eqnarray*}
\end{proof}

\subsection{Prym-Tyurin varieties} Recall that an abelian subvariety $P$ of $JX$ is called a {\it Prym-Tyurin variety of exponent}
$q$ in $JX$, if the restriction of the canonical principal polarization of $JX$ to $P$ is the $q$-fold of a principal polarization of $P$.
The main result of the paper is the following theorem.

\begin{thm} \label{thm4.8}
Let $\cW_1, \ldots, \cW_r$
denote nontrivial pairwise non-isomorphic rational irreducible representations of the group $G$
with associated complex irreducible representations $V_1, \ldots, V_r$ satisfying Hypothesis \ref{hyp}.
Suppose that the action of the finite group $G$ has geometric signature $[0;(C_1,m_1), \ldots , (C_t,m_t)]$
satisfying
\begin{equation}  \label{eqn4.3}
\sum_{j=1}^t m_j \left( q[L:\QQ]\sum_{i=1}^r  (\dim V_i - \dim V_i^{G_j})  -
([G:H] -|H \backslash G/G_j|) \right) = 0.
\end{equation}
where $G_j$ is of class $C_j$ and $q = \frac{|G|}{b \cdot n}$.  Then $P_{\overline \cD}$ is a Prym-Tyurin
variety of exponent $q$ in $JX$.
\end{thm}

Note that the theorem gives a method to construct
Prym-Tyurin varieties.

\begin{proof}
Since $g(Y) = 0$, we have $P(X/Y) = JX$. The correspondence $\cK_X$ on the curve $X$ is integral, symmetric and fixed-point free according to Proposition \ref{prop4.6}. Moreover
the associated endomorphism $\kappa_X$ of $JX$ satisfies equation \eqref{eq4.1}. Hence the assertion is a consequence
of Kanev's criterion (see e.g. \cite[Theorem 12.9.1]{bl}).
\end{proof}

\subsection{Relation to Kanev's construction}
In this subsection we compare our construction of the Prym variety $P_{\overline \cD}$ to Kanev's original construction. Let us first recall the original
construction.

Let $\cW$ be an absolutely rational irreducible representation of the group $G$. Suppose $G$ acts on a lattice
$\Lambda$ of maximal rank in $\cW$. Moreover let $w \in \cW$ be a weight; that is a nonzero vector satisfying $gw - w \in \Lambda$ for all $g \in G$.
Then there exists a uniquely determined negative definite $G$-invariant symmetric bilinear form $( \;|\;)$ on $\cW$ such that
(i) $(w | \Lambda ) \subset \ZZ$ and (ii) any negative definite $G$-invariant form with (i) is an integer multiple of $( \;|\;)$.

Now let $\pi: Z \ra \PP^1$ be a simply ramified Galois covering with Galois group $G$. Consider the subgroup  $H := \Stab_G(w)$,
define $X := Z/H$, and let $\varphi_H: X \ra \PP^1$ denote the canonical map.
Then for any pair $(\Lambda, w)$ as above, Kanev defines a symmetric effective correspondence
on the curve $X$, which induces
an abelian subvariety $P_K = P_K(\Lambda,w)$ of $JX$. We call it the {\it Prym variety} in $JX$ associated to the pair $(\Lambda,w)$.
It is shown in \cite{lr} that the abelian subvariety $P_K$ of $JX$ is given by the image of the element $s_w \in \End_{\QQ}(JX)$
associated to the correspondence $S_w(x)$ on $X$, defined for every $x \in X$ with
preimage $z \in Z$ by (see \cite[Proposition 2.2]{lr} and its proof)
\begin{equation} \label{eq4.3}
S_w(x) = |H|^2 \sum_{i=1}^d \sum_{j=1}^{n_i} (w|g_{ij}w)\pi_H(g_{ij}(z)).
\end{equation}
Here $\{g_{ij} \; | \; i = 1, \ldots,d;\; j = 1 , \ldots, n_i \}$ is a set of representatives of the left and right cosets of $H$ in $G$.

%Fix a point $\xi_0$ of $\PP^1$ not in the
%branch locus of $\varphi_H$. Choosing an element in the fibre $\varphi_H^{-1}(\xi_0)$ induces a bijection
%$$
%\{g_1w, \ldots, g_dw\} \stackrel{\sim}{\ra} \varphi_H^{-1}(\xi_0)
%$$
%which is $G$-equivariant according to the definitions. Here $\{g_1 = 1, g_2, \ldots, g_d\}$ is a set of presentatives
%simultaneously for the left and right cosets of $H$ in $G$. In the sequel we identify both sets, i.e. we label the elements of
%$\varphi_H^{-1}(\xi_0)$ by $g_1w, \ldots g_dw$.

Now let the notation be as in Section \ref{S:4}. Consider the special case
$r = 1$, $V = \cW \otimes \CC$ with an absolutely rational irreducible representation $\cW$ and finally $H$ a subgroup of $G$
satisfying $\dim V^H = 1$.
Clearly $\dim \cW^H =1$. For any
nonzero vector $w \in \cW^H$ the subgroup of $\cW$ generated by the elements $gw -w$ for all $g \in G$ is a $G$-invariant lattice
$\Lambda$ of maximal rank in $\cW$ for which $w$ is a weight. Hence the Prym variety associated to the pair $(\Lambda,w)$
is well defined. With these assumptions we have

\begin{prop} \label{prop4.9}
The Prym variety associated to the pair $(\Lambda,w)$ coincides with the Prym variety $P_{\overline \cD}$ as defined in Section \ref{SS:prym}:
$$
P_K(\Lambda,w) = P_{\overline \cD}.
$$
\end{prop}

\begin{proof}
Note first that $P_{\overline \cD}$ does not depend on the choice of the bilinear form $(\cdot, \cdot)$
used in Section \ref{S:2}. The statement of
Proposition \ref{prop2.5} is independent of the choice of the form $(\cdot , \cdot)$. Hence we may choose the bilinear form in such a way that
$$
(w,w) = \frac{|H|}{b},
$$
with $b = \gcd\{a_1 -a_i \;|\; 2 \leq i \leq d \}$ as in \eqref{eqn3.5}.

Then we have, according to Proposition \ref{prop2.4} and equations \eqref{eqn2.3} and \eqref{eq2.6},
\begin{equation} \label{eq4.4}
(w,g_{ij}w) = \frac{(w,w)}{|H|} \sum_{h \in H} \chi_V(hg_{ij}^{-1})
= \frac{(w,w)}{|H|} \sum_{h \in H} \chi_V(hg_{i1}^{-1})  = \frac{a_i}{b}
\end{equation}
and thus,  for all $i$ and $j$,
$$
(w,g_{ij}w - w) = \frac{a_i - a_1}{b}.
$$
This implies that Kanev's distinguished form $( \cdot | \cdot)$ is just the negative of our form $( \cdot , \cdot )$.

Now the Prym variety $P_{\overline \cD}$ is defined as the image of the element $\delta_{\overline \cD} \in \End_{\QQ}(JX)$
associated to the correspondence
${\overline \cD}$, which according to equation \eqref{e3.5} is given by
\begin{eqnarray*}
{\overline \cD}_{H,\cW}(x) & = &\sum_{i=1}^d a_{i} \sum_{j=1}^{n_i} \pi_H g_{ij}(z)\\
& = &  \sum_{i=1}^d b (w, g_{ij}w) \sum_{j=1}^{n_i} \pi_H g_{ij}(z) \quad  (\mbox{by} \eqref{eq4.4})\\
& = & -b \sum_{i=1}^d  \sum_{j=1}^{n_i} (w| g_{ij}w) \pi_H g_{ij}(z).
\end{eqnarray*}
Comparison with \eqref{eq4.3} implies that $s_w$ and $\delta_{\overline \cD}$ are rational multiples of each other, and
thus have the same image in $JX$.
\end{proof}

\begin{rem}
(a): Proposition \ref{prop4.9} means that Kanev's and our constructions have nonzero intersection. Kanev's construction
is more general in the sense that he does not assume that $\dim V^{\Stab(w)} = 1$. Our construction is more general in the sense
that we do not have to assume $\cW$ absolutely irreducible. Moreover, $\cW$ can be a tuple of representations, which actually gives
interesting examples, as we will see in Section \ref{SS:oe}.\\
(b): Our construction of abelian subvarieties of $JX$ can be carried out without the hypotheses ``$\dim V^H = 1$
and $H$ maximal with this property'' and ``$Z/G = \PP^1$'' (and this is valid in fact for any
rational irreducible representation $\cW$) by just considering a primitive idempotent $f$
in $Q[G]e_{\cW}$ that is invariant under multiplication by $p_H$ on both
sides; it can be proven that such idempotents exist if $\langle \rho_H , V \rangle
\neq 0$. Then $P = \im(f)$ is an abelian subvariety of $JX$.
Hence our construction also has points in common with the one by Merindol \cite{M},
again generalizing to non absolutely irreducible representations.

The purpose of our additional assumptions is to allow the study of the
polarization of $P$.
\end{rem}

\section{Examples}\label{S:5}

\subsection{Kanev's Examples}
In \cite{k} Kanev constructed families of Prym-Tyurin varieties for the Weyl groups $G = W(R)$ for the root systems $R$
of type $A_n, D_n, E_6$ and $E_7$ using the construction outlined at the beginning of Section 4.6. In these cases $\cW$
is the root representation of the group $G$,  $\Lambda$ is the root lattice and $w$ a minuscule fundamental weight.
It is easy to check that in all these cases the subgroup $H = \Stab_G(w)$ satisfies $\dim \cW^H = 1$ and is maximal
with this property. So, according to Proposition 4.9 the corresponding Prym-Tyurin varieties can be also constructed
by our method.

\subsection{Additional Prym-Tyurin varieties for Weyl groups}
In order to construct his Prym-Tyurin varieties, Kanev \cite{k} starts with a covering $\varphi_H: X \ra \PP^1$
with simplest ramification. In our terminology this means that the geometric signature of the covering $\pi: Z \ra Y$
only contains the conjugacy class of reflections in the Weyl group. Using equation \eqref{eqn4.3} one can work out exactly the types of
conjugacy classes to which Theorem \ref{thm4.8} applies. We call these classes {\it admissible} conjugacy classes for the triple
$(G,H,V)$, where $G = W(R)$, $V$ is the root representation and $H$ a subgroup of $G$ satisfying Hypothesis \ref{hyp}.
Using coverings with these more general geometric signatures we get other families of Prym-Tyurin varieties. In fact they
are of different dimensions. However certainly these varieties can also be constructed generalizing Kanev's approach.\\

One idea to construct new Prym-Tyurin varieties is as follows:
For each root system $R$ of type $A_n$, $D_n$, $E_6$ and $E_7$,
denote by $G$ its Weyl group. In each of these cases $G$ has a
unique subgroup of index $2$, which will be denoted by
$\widetilde{G}$; in fact, if we choose the basis $\{ \alpha_i\}$ for
the corresponding root system as in Bourbaki \cite{b}, then $G$ is generated
by the reflections $r_{\alpha_i}$, and $\widetilde{G}$ is generated
by the elements of order three in $G$ that are products of two such
reflections $s_{ij} = r_{\alpha_i} r_{\alpha_j}$ for $\alpha_i$ and
$\alpha_j$ connected by an edge in the corresponding Dynkin diagram. Let us denote
by $C_2$ the class of a subgroup $G_2$ generated by one of the $r_{\alpha_i}$, and by
$C_3$ the class of a subgroup $G_3$ generated by one of the $s_{ij}$.

\begin{prop}
Assume $Z$ is a curve with $\widetilde{G}$-action and geometric signature $[0;(C_3,m)]$.
Let $V$ denote the restriction to $\widetilde{G}$ of the root
representation of $G$.

Then the triples $(\widetilde{G}, H , V)$, with $H$ as in the
following table, satisfy Hypothesis \ref{hyp}, and the corresponding
abelian subvarieties of $J(X), \; X = Z/H$, are Prym-Tyurin varieties of exponent $q$ as
in the table.

\begin{center}
\begin{tabular}{|c|c|c|}
% \hline
$\widetilde{G}$  & $H$     & $q$       \\
\hline
%   & &      \\ \hline
   & &     \\
$\widetilde{WA_n} = \Alt(n+1)$  & $(S_k\times S_{n+1-k}) \cap \Alt(n+1)$ , $1 \leq k \leq n$   & $\frac{(n-1)!}{(k-1)!(n-k)!}$   \\
%    &     &      \\
   &    &   \\
$\widetilde{WD_n}$ &
$\left( (\mathbb{Z}/2\mathbb{Z})^{n-2} \rtimes S_{n-1}\right) \cap \widetilde{WD_n}$ & $2$ \\
%   &     &                             \\
    &    &                           \\
$\widetilde{WD_n}$ & $S_n \cap \widetilde{WD_n}$  & $2^{n-3}$ \\
%  &   &     \\
  &  &         \\
$\widetilde{WE_6}$  & $\widetilde{WD_5}$ & $6$\\
% &   &       \\
  &  &      \\
$\widetilde{WE_7}$ & $\widetilde{WE_6}$ & $12$\\
\hline
%  &  &       \\ \hline
\end{tabular}
\end{center}
\end{prop}

\begin{proof} The proof is a straightforward application of Theorem \ref{thm4.8}.
We use for it the computer program mentioned in the introduction. We omit the details.
\end{proof}

Note that we have obtained new families of Prym Tyurin varieties, with the same
exponent as in the examples of Section 5.1, and that the representation $V$ is in each
case absolutely irreducible.
By considering other representations of the same groups, we can find
examples with
other exponents, and also cases where the restriction is not
absolutely irreducible. Some examples are given below.
\\

The group $WD_5$ has a unique representation of degree six; it is
absolutely irreducible, and there is no Prym-Tyurin variety associated to this
representation that may be constructed using our methods. However,
when restricted to $\widetilde{WD_5}$  it becomes a rational
representation $\mathcal{W}$ with associated complex representation
$V$ of degree three and $[K_V: \mathbb{Q}] = 2$. Furthermore, there
exists a (unique conjugacy class of) subgroup $H$ of order $80$ of
$\widetilde{WD_5}$ such that $(\widetilde{WD_5}, H , \mathcal{W})$
satisfies Hypothesis \ref{hyp}. Again we omit the details of the proof of the following proposition.

\begin{prop}
Consider the triple $(\widetilde{WD_5}, H , \mathcal{W})$, then the conjugacy classes
admissible for it are the
two conjugacy classes $C_6,C_6'$ of cyclic subgroups of  order $6$ and the
unique conjugacy class $C_3$ of order $3$, each of them generates the group $\widetilde{WD_5}$.
Moreover, if $m_3,m_6,m_6'$ are the number of branch points of type $C_3,C_6,C_6'$ respectively, then
the associated Prym-Tyurin variety has exponent $2$ and dimension
$2(m_3+m_6+m_6')-6$.

Similarly, for the irreducible rational representation $\mathcal{W}_4$ of dimension $4$ of
$\widetilde{WD_5}$, there exists a subgroup $H_2$ of order $96$ such
that $(\widetilde{WD_5}, H_2 , \mathcal{W}_4)$ satisfies Hypothesis
\ref{hyp}. The same conjugacy classes as above are admissible for this triple, and the associated
Prym-Tyurin variety has exponent $3$ and dimension $m_3+m_6+m_6'-4$.
exponent $3$ and dimension $m_3+m_6+m_6'-4$.
\end{prop}

\subsection{Other examples}\label{SS:oe}

As mentioned in the introduction, we obtain many new families of
Prym-Tyurin varieties, if we consider the case $r>1$ with the notation of Hypothesis \ref{hyp}.
We only mention two series. The first one shows how to construct Salomon's examples via Theorem \ref{thm4.8}.\\

Consider the triple $(G = S_n \times S_n, H , \mathcal{W} =(V \boxtimes
V_0 , V_0 \boxtimes V))$, where $H$ is a subgroup of index $n^2$, $V$
is the root representation of $WA_{n-1} = S_n$, and $V_0$ is the trivial
representation of $S_n$. This triple satisfies the Hypothesis \ref{hyp}.

By choosing the geometric signature of type
$[0;(C_1,m_1),(C_2,m_2)]$ where $C_1$ is the class of $G_1 = \langle ((1\,2),1_{S_n}) \rangle$ and
$C_2$ the class of
$G_2 = \langle (1_{S_n},(1\,2))\rangle$, we reobtain Salomon's examples of
Prym-Tyurin varieties of exponent $n$. If we choose a different signature, and there
are many possibilities (for instance another admissible group is
$G_3 = \langle ((1\,2),(1\,2)) \rangle$, we obtain new families of Prym-Tyurin varieties
of exponent $n$ for each natural $n$.
Examples may also be constructed with non absolutely irreducible
representations.
\\

For the second series of examples let $D_p$ denote the dihedral group of order $2p$ with an odd prime $p$.
It has $\frac{p-1}{2}$ complex
irreducible representations of degree two, each with $L_V = K_V =
\mathbb{Q}[w_p+w_p^{-1}]$, where $w_p$ denotes a $p$-th root of
unity. They are all associated to the same rational irreducible
representation. Call any one of them $V$.

Now consider the group
$$
G = D_p \times D_p \times \mathbb{Z}/2\mathbb{Z}= \;\langle x_1, y_1
: x_1^p, y_1^2, (x_1y_1)^2 \rangle \times \langle x_2, y_2 : x_2^p,
y_2^2, (x_2y_2)^2 \rangle \, \times \langle z : z^2\rangle ,
$$
the subgroup
$$
H = \langle y_1 , y_2 , z\rangle ,
$$
and the pair of rational irreducible representations
$\mathcal{W}_j$, $j=1, 2$, of $G$  associated to $V\boxtimes V_0
\boxtimes V_1$ and $V_0 \boxtimes V \boxtimes V_1$ respectively, where $V$
is as above, and $V_0$ and $V_1$ denote the trivial representations
for $D_p$ and  $\mathbb{Z}/2\mathbb{Z}$ respectively.
Then one checks that the triple
$(G,H,\mathcal{W}=(\mathcal{W}_1,\mathcal{W}_2))$ satisfies the
conditions of Hypothesis \ref{hyp}.
Furthermore, the classes of the subgroups $G_1  = \langle
y_1 \rangle$,  $G_2  = \langle y_2 \rangle$, $G_3  = \langle y_1\,
z\rangle$,   $G_4  = \langle y_2\, z\rangle$, $G_5  = \langle x_1
\rangle$,   $G_6  = \langle x_2 \rangle$, $G_7  = \langle x_1 \,
z\rangle$, $G_8  = \langle x_2 \, z\rangle$, and $G_9  = \langle z
\rangle$ are admissible for the triple $(G,H,\mathcal{W})$.
Then one can apply Theorem \ref{thm4.8} to deduce
the following proposition, where for the details of the proof we refer to \cite{clrr}.

\begin{prop}
Let the notation be as above and let $Z$ be a curve with $G$-action and geometric
signature $[0;(G_1,m_1) , \ldots ,$
$ (G_9,m_9)]$ with $m_i$ even for $1 \leq i \leq 8$.
Then, if $\ell_1 := m_1 + m_3 + 2(m_5+m_7)$ and $\ell_2 := m_2 + m_4 + 2(m_6 + m_8)$ are both $ \geq 6$,
the associated abelian subvarieties of $JX$, with $X = Z/H$ as usual, are
Prym-Tyurin varieties of exponent $p$ and dimension $\frac{p-1}{4}
(\ell_1+\ell_2-8)$.
\end{prop}

\end{document}